\documentclass{amsart}
\date{}
\usepackage{amsmath}
\usepackage{amscd}
\usepackage{amssymb}
\usepackage{amsfonts}
\usepackage{latexsym}
\usepackage{verbatim}
\usepackage{epsfig}

\numberwithin{equation}{section}

\newtheorem{problem}[equation]{Problem}
\newtheorem{proposition}[equation]{Proposition}
\newtheorem{theorem}[equation]{Theorem}

\newtheorem{lem}[equation]{Lemma}
\newtheorem{prop}[equation]{Proposition}

\theoremstyle{remark}

\newtheorem{remark}[equation]{Remark}
\newtheorem{definition}[equation]{Definition}

\newcommand{\ignore}[1]{}

\newcommand{\eps}{\epsilon}

\newcommand{\N}{\mathbb N}
\newcommand{\QED}{$\Box$}
\newcommand{\R}{\mathbb R}

\newcommand{\Z}{\mathbb Z}

\newcommand{\annotation}[1]{\marginpar{\tiny #1}}

\begin{document}

\title
[Potentials and almost isomorphism ] {Good potentials for almost
isomorphism of countable state Markov shifts}

\begin{abstract}
Almost isomorphism is an equivalence
relation  on countable state Markov shifts
which provides a  strong version of Borel conjugacy;
 still, for mixing SPR shifts,
entropy is a complete invariant of almost isomorphism  \cite{BBG}.
In this paper, we establish a class of potentials on countable
state Markov shifts whose thermodynamic
formalism is respected by almost isomorphism.
\end{abstract}

\author{Mike Boyle}
\address{Mike Boyle\\Department of Mathematics\\
        University of Maryland\\
        College Park, MD 20742-4015\\
        U.S.A.}
\email{mmb@math.umd.edu} \urladdr{www.math.umd.edu/$\sim$mmb}
\author{Jerome Buzzi}
\address{Jerome Buzzi\\Centre de Math\'ematiques\\
Ecole polytechnique\\91128 Palaiseau Cedex\\France}
\email{buzzi@math.polytechnique.fr} \urladdr{www.jeromebuzzi.com}
\author{Ricardo G\'omez}
\address{Ricardo G\'omez\\ Instituto de Matematicas Area de la Investigacion Cientifica\\
Circuito Exterior\\ Ciudad Universitaria\\ DF 04510\\ Mexico }
\email{rgomez@math.unam.mx}
\urladdr{www.math.unam.mx/$\sim$rgomez}

\thanks{This work was
partly supported by NSF Grant 0400493 (Boyle). }




\maketitle

\tableofcontents

\section{Introduction}

The classical subshifts of finite type \cite{BowenLNM,Ki,LM} admit the
following natural generalization. A {\em countable state Markov
shift} is $(X,S)$ where $X\subset V^\Z$ for some countable (maybe
finite) set $V$ and for some $E\subset V^2$:
 $$
    X = \{ x\in V^\Z:\forall n\in\Z\; (x_n,x_{n+1})\in E\}
 $$
and $S:X\to X$ defined by
 $$
    S((x_n)_{n\in\Z})=(x_{n+1})_{n\in\Z}.
 $$
In other words, $X$ is the set of bi-infinite paths on the
directed graph $(V,E)$ together with the left-shift $S$. The
one-sided version $(X_+,S_+)$ is given by $X_+=\{x_0x_1\dots\in
V^\N:x\in X\}$ and $S_+((x_n)_{n\in\N})=(x_{n+1})_{n\in\N}$.
(We shall sometimes use the same letter for both the map and the space.)
A countable state Markov shift is
{\it transitive}, or
{\it irreducible},  if for any vertices $u,v$ in the underlying directed
graph there is a path from $u$ to $v$; and it is
{\it mixing} if  for given $(u,v)\in V^2$, for all but finitely many
$k$, this path may be chosen to have length $k$.
The {\it topological entropy}, or {\em Gurevi\v{c} entropy} \cite{G0},
 of $(X,S)$
is $h(S)=\sup_\mu h(S,\mu)$, the supremum of the entropies of
$S$-invariant probability measures (see \cite{Walters-book} for
background on entropy).

STANDING CONVENTION: for the rest of this
paper,  unless there is an explicit qualification,
 ``Markov shift'' means
 an irreducible  countable (maybe finite) state Markov shift of
finite topological entropy.

A Markov shift $(X,S)$ is {\em strongly positive recurrent} or
SPR \cite{Vere-Jones2,G1,UF,salama2,ruette,QFT} if it admits an
invariant probability measure $\mu$ with $h(S,\mu)=h(S)$ and
\[
\limsup_{n\to\infty} (1/n)\log\mu\left(X\setminus \bigcup_{k=0}^n
S^{-k}U\right)<0
\]
 for all nonempty open $U\subset X$.
In the terminology of \cite{GS}, the
SPR shifts  are the positive recurrent
Markov shifts whose defining directed
graphs have adjacency matrices which are
{\it stable positive}.
There are other characterizations of SPR Markov shifts (assembled
in \cite[Prop. 2.3]{BBG}). The SPR shifts are the natural large
class of Markov shifts
which still retain certain key properties of shifts of finite type.

{\em Almost isomorphism} of Markov shifts yields a strong version
of Borel conjugacy (recalled in Section \ref{sec:statement}). In
particular, almost isomorphism of Markov shifts  induces a
bijection between their sets of ergodic fully supported measures,
simultaneously defining isomorphisms of all the corresponding
measurable systems. This gives an identification of the measures
of maximal entropy, when they exist, by a map for which the coding
time (in the SPR case) has an exponential tail.
Our main result in \cite{BBG} shows that for mixing SPR Markov
shifts, topological entropy is a complete invariant of almost
isomorphism.

Thermodynamic formalism generalizes the notion of measure of
maximal entropy to equilibrium measure (\ref{defneqmeas}) of  a
function (potential).  An equilibrium measure for a
``nice'' potential  should be fully supported and ergodic
(remember that we assume our systems to be irreducible). Almost
isomorphism respects the class of such measures, so this begs the
question, is there a class of reasonably nice potentials which
together with their equilibrium measures are respected by almost
isomorphism?
 The usual classes of
nice potentials (locally H\"older continuous or with summable
variations) are not preserved by almost isomorphism.

In this paper we introduce a new class of potentials called {\em
relatively regular} which, on the one hand, are nice enough to
guarantee existence of a unique equilibrium measure which is fully
supported and, on the other hand, are invariant under almost
isomorphism. This class contains the
 positive recurrent potentials of summable
variation, but necessarily contains certain less regular
potentials as well (though it does not include all potentials
which even smooth functions can generate by coding when there is
non-uniform expansion).
 Our main result is that almost
isomorphisms leave globally invariant the thermodynamic formalism
of the bounded relatively regular functions.

 \ignore{Markov shifts with infinitely many states arise
naturally in the coding of nonuniformly hyperbolic smooth or
piecewise smooth systems (see e.g. \cite{H3,K2,Young};  for SPR
especially, see e.g. \cite[Sec.7]{BBG} and \cite{AFF,SIM,QFT}). In
coding applications, one may be interested not only in the measure
of maximal entropy, but also in equilibrium measures of classes of
functions (potentials) \cite{BowenLNM,  Young}. Consequently, the
aim of this paper is to establish a reasonably large class  of
potentials with unique equilibrium measures, such that the
potentials and their equilibrium measures are respected by almost
isomorphism. The class we choose  is the class of bounded
``relatively regular'' potentials (\ref{goodclass}). This class
contains the bounded positive recurrent potentials of summable
variation, but it must then be enlarged to contain certain less
regular potentials, because properties like
  summable variation (or H\"{o}lder)
 are not respected by almost isomorphism. We choose the class small
enough that we can still prove existence and uniqueness of
equilibrium measures.}

In Sec \ref{sec:statement}, we provide definitions and state
our main result, Theorem \ref{boundedtheorem}: the thermodynamic
formalism of the bounded relatively regular functions is
respected by almost isomorphism.
 The uniqueness of their equilibrium measures is
established (in greater generality) in Section
\ref{sec:uniqueness}; the existence is established (in greater
generality) in \ref{sec:existence}; and the correspondence under
almost isomorphism is established in Section
\ref{sec:correspondence}. \ignore{ In Section \ref{sec:example} we
interpret the results in the setting
  of certain interval maps, and ask a question.}


\section{Statement of result}\label{sec:statement}

Given a Borel map $S:X\to X$,  a real-valued Borel measurable
function $f$ from $X$ to $\mathbb R$, and an $S$-invariant
Borel probability $\mu$, we define
\begin{align*}
P(S,f,\mu ) &= h_{\mu }(S) + \int f\ d\mu \ ,\ \ \ \ \textnormal{ and} \\
P(S,f)      &= \sup_{\mu} P(S,f,\mu ) \ .
\end{align*}
We take $P(S,f)$ as our definition of the {\it pressure} of $(S,f)$.
 $P(S,f,\mu )$ might not be defined for some $\mu$; the supremum
defining $P(S,f)$ is taken over the well defined $P(S,f,\mu )$.
\begin{definition}\label{defneqmeas}
An {\it equilibrium measure} of $f$ (also called an equilibrium state)
is an $S$-invariant Borel probability $\mu$ such that
$P(S,f,\mu ) =P(S,f)$.
\end{definition}
For background see \cite{Walters-book,GS,FFY,MU1,Sarig0}.

Next we recall the definition of almost isomorphism from
\cite{BBG}. A map $\varphi :S\to T$ between Markov shifts is a
{\it one-block code} if there is a function $\Phi $ from the
symbol set of $S$ into the symbol set of $T$ such that $(\varphi
x)_n = \Phi  (x_n)$, for all $x$ and $n$  (note that $\varphi
x=\varphi(x)$). A $T$-word $W$ (of length $|W|$) is a {\it magic
word} for such a map $\varphi $ if the following hold. We denote
by $x[a,b]$ the restriction of the sequence $x$ to the indices
$i=a,a+1,\dots,b$.
\begin{enumerate}
\item \label{mw1} If
 $y\in T$ and
$\{n\in \mathbb Z: y[n,n+|W|-1]=W\}$ is unbounded above and
unbounded below, then $y$ has a preimage under $\varphi $. \item
\label{mw2} There is an integer  $I$ such that whenever $C$ is a
$T$-word and two points $x$ and $x'$ of $S$ satisfy $(\varphi
x)[0,2|W|+|C|-1]=WCW= (\varphi  x')[0,2|W|+|C|-1]$, then $x[I,
I+|W|+|C|-1]= x'[I, I+|W|+|C|-1]$.
\end{enumerate}
(In constructions, the integer $I$ of the last condition can
generally be chosen to be zero.) It follows from (\ref{mw2}) that
the preimage in (\ref{mw1}) is unique.

\begin{definition}\label{def-ai}
 Markov shifts $S$ and $T$  are
 {\em almost isomorphic} if there exist a
Markov shift $R$ and one-block codes $ R\to
S$, $R\to T$ each of which is injective with a
  magic
word. Such a pair of maps is an
{\it almost isomorphism} of $S$ and $T$.
\end{definition}

Our interest in almost isomorphism is largely explained by
the following proposition, copied from \cite[Proposition 3.4]{BBG}.
(We  use ``vertex shifts'' in this paper rather than the
``edge shifts'' of \cite{BBG}, but this is only a matter of
notation.)

\begin{proposition}\label{magicprop}
Suppose $S$ and $T$ are almost isomorphic Markov shifts. Then
$h(S)=h(T)$, and there are Borel subsets $K$ and $K'$ of $S$ and
$T$,  collections of invariant probability measures
$\mathcal M(K),\mathcal M(K')$ on $K$ and $K'$ and a
shift-commuting Borel-measurable bijection $\gamma :K\to K'$, such
that the following hold.
\begin{enumerate}
\item \label{ai1} $K$ and $K'$ are residual subsets of $S$ and $T$
(contain dense $G_{\delta}$ sets). \item \label{ai2} The map
$\gamma$ induces a bijection
 $\mathcal M(K)\to \mathcal M(K')$ ($\mu \mapsto \mu '$, say)
such that for each such pair $\mu,\mu'$
 the map $\gamma$ induces an
 isomorphism $\gamma :(S,\mu )\to (T, \mu ')$,
which is a magic word isomorphism when $\mu$ and $\mu'$ have full
support.
 \item
\label{ai3} $\mathcal M(K)$ and $\mathcal M(K')$
 contain all ergodic shift-invariant Borel
probabilities on $S$ and $T$ with full support, and these
correspond under $\gamma$. \item \label{ai5} If $S$ is SPR, then
so is $T$, and $\gamma $ is an  entropy-conjugacy from $S$ to $T$.
\end{enumerate}
\end{proposition}

In the proposition, ``full support'' means that the
measure is  nonvanishing on all nonempty open sets, and
 ``entropy-conjugacy''
means that there
exists $\epsilon >0$ such that the sets $\mathcal M(K)$,
$\mathcal M(K')$ contain all invariant ergodic Borel
probabilities with entropy greater than $h(S)-\epsilon$.
We denote by
$\mathcal M^{\text{erg}}_{\text{supp}}(S)$ the set of ergodic
fully supported $S$-invariant Borel probabilities.

Given a Markov shift $S$ and a nonempty collection $\mathcal W$ of
$S$-words (words appearing in points of $S$), let $S_{\mathcal W}$
denote the set of all points $x$ in $S$ which see $\mathcal W$
words infinitely often in the past and future\footnote{Dealing
with ergodic invariant probability measures, we could equally
require that \emph{some} $\mathcal W$-word is seen infinitely
often in the past and future.}. (When $\mathcal W=\{W\}$ we may
use the notation $S_W$.) If $ W $ is a magic word for the map
$R\to S$ of Definition \ref{def-ai}, then the magic word
isomorphism  of Definition \ref{def-ai} induces an obvious
shift-commuting Borel injection $\gamma$ from $S_{ W}$ into $T$.
The word $W$ can be chosen so that $\gamma$ has an inverse
similarly defined on a subset of some $T_{\mathcal W'}$ in $T$,
and then $S_{W}$ can be used for the set $K$ in Proposition
\ref{magicprop}. For lighter notation, we will generally use the
same symbol $\gamma$ for a map defined in this way and for its
induced maps on functions and measures.

Properly speaking, then, almost isomorphisms relate not functions in
themselves but rather equivalence classes, as follows.
Given Borel functions $f$,$g$
defined on subsets of $S$, we say $f$ and $g$ are {\it somewhere
equivalent} ($f=g$ s.e.)  if there exists some $S_{\mathcal W}$
such that $f$ and $g$ are defined and equal on all of $S_{\mathcal
W}$.  When we say below that two functions correspond mod s.e.
under some almost isomorphism $R\to S$, $R\to T$, we mean that
some representatives $f,g$ of their s.e. equivalence classes
correspond, i.e., satisfy $g=f\circ\gamma$ s.e. where $\gamma$ is defined
as above.

In considering the behavior of potentials and equilibrium measures
under almost isomorphism, we are interested in classes of
functions invariant under almost isomorphism; in particular, these
functions should be described by properties which persist under
restriction to smaller sets $S_W$. Note,
 if a measure does
not have full support, then it will assign measure zero to some
$S_W$. Conversely, if $f =g$ s.e., then $f=g$ $\nu$-a.e. for all
$\nu$ in $\mathcal M^{\text{erg}}_{\text{supp}}(S)$. Therefore we
are interested in functions all of whose equilibrium measures have
full support.

 First we define some regularity conditions related
to full support (and uniqueness) of equilibrium measures. Let $S$
be a Markov shift and let $\mathcal W$ denote a set of nonempty
$S$-words.
\begin{definition}\label{psummdef}
A function $f:S\to \mathbb R$ has
{\it eventually  $p$-summable variations relative to $\mathcal W$}
if there exists a sequence $\omega_1\geq\omega_2\geq\dots$ with
 $$
   \sum_{n\geq1} n^p \omega_n < \infty
$$
such that for all $x,y$ in $S$ and integers $m,n\geq0$, we have
$|f(x)-f(y)|\leq \omega_n+\omega_m$ whenever
\begin{enumerate}
 \item\label{condi} $x[-m,n]=y[-m,n]$,
 \item\label{condii} $x[-m,-1]$ begins with
a word from $\mathcal W$, and
 \item\label{condiii} $x[1,n]$ ends with a word from
$\mathcal W$.
\end{enumerate}
\end{definition}
We remark that conditions (\ref{condii}) and (\ref{condiii}) imply
that $m,n\geq1$. On the other hand, if we removed these two
conditions and set $\omega_n=\kappa^n$ for some $0<\kappa<1$, we
would be describing H\"{o}lder continuity with respect to a
suitable metric -- a metric ``of the first type'' (1.2) in
\cite{GS}.

Given $p\geq 0$, we define $\mathcal E_p(S,\mathcal
W)$ as the set of functions with eventually $p$-summable
variations relative to $\mathcal W$, and set
\[
\mathcal E_p(S) =\bigcup_{\mathcal W}\mathcal
E_p(S,\mathcal W)\ \ .
\]
Like  the functions with summable
variations in \cite{Sarig0}, the functions in $\mathcal
E_p(S)$ are not necessarily bounded.
We let $\mathcal E_{p+}(S)$ denote the subset of
$\mathcal E_p(S) $ consisting of functions
 depending only on future coordinates, i.e.,
$f(x)=f(x_0,x_1,\dots )$. For such $f$, the conditions of
Definition \ref{psummdef} simplify (because any word can be
continued to the left to a word in $\mathcal W$): if
$x[0,n]=y[0,n]$, with $x[1,n]$ ending with a word from $\mathcal W$,
then
 $|f(x)-f(y)|\leq
\omega_n$.

The regularity of the potential is not sufficient to guarantee
existence of an equilibrium measure (even the zero potential may
fail to have an equilibrium measure \cite{G0}). For existence, we
will use a notion from \cite{Sarig0} (generalized from \cite{G0}).
The {\it local partition function} of $(S,f)$ at some cylinder
$[W]$ associates to each positive integer $n$ the sum
 $$
    Z_n(S,f,W)=\sum_{\begin{matrix} S^nx=x\\x[0,|W|-1]=W\end{matrix}}
      \exp \left( f(x)+f(Sx)+\dots+f(S^{n-1}x) \right)\ .
 $$
The pair $(S,f)$ is said to be {\sl positively recurrent} if for
some (or equivalently, every) nonempty $S$-word $W$
there exists an
integer $n_0=n_0(W)$ such that
the  sequence
 $$
     \Big(Z_n(S,f,W)\exp (-nP(S,f) )\Big)_{n\geq n_0}
 $$
is bounded away from $0$ and $\infty$.

\begin{definition}\label{goodclass}
$ \mathcal E(S)$
denotes the set of {\it relatively regular} functions on $S$.
These are the
positive recurrent Borel measurable
functions $f$ from $S$ to $\mathbb R\cup\{-\infty\}$ such that
$P(S,f)<\infty$ and there exists a nonempty family
 $\mathcal W$ of $S$-words (which may depend on $f$) such that
 $f\in \mathcal E_1(S,\mathcal W)\cup
\mathcal E_{0+}(S,\mathcal W )$ and every ergodic
equilibrium measure of $f$ assigns measure $1$ to $S_{\mathcal
W}$. The set of bounded relatively regular functions is denoted
$ \mathcal E^{\flat }(S)$.
\end{definition}
Let us point out that in the special case that $\mathcal
W=\mathcal A(S)$ (the alphabet of $S$), the  condition in
Definition \ref{goodclass} that $S_\mathcal W$ has full measure
follows from the others (see \cite{Sarig0}).

 The following theorem, our main result, shows that
the bounded relatively regular functions comprise
one good class of potentials for almost
isomorphism.

\ignore{\annotation{I suppressed the confessional remark as
suggested ... Some references (by Jerome, Hofbauer, Keller, Young)
were no longer cited in the text, so I removed those from the
bibliography, except for the QFT paper, for which I put a citation
in with the SPR list.}}

\begin{theorem}\label{boundedtheorem}
 Suppose $R\to S$, $R \to T$ is an almost isomorphism of
Markov shifts $S$ and $T$. This almost isomorphism induces a
mod-s.e. bijection
 $\ \gamma: \mathcal E^{\flat }(S) \to
\mathcal E^{\flat }(T) $.
For every  $f\in  \mathcal E^{\flat }(S)$, there is
 a unique
equilibrium measure, $\mu_f$; this measure $\mu_f$ has full
support; and the almost isomorphism induces a
correspondence $\gamma : \mu_f \to \mu_{(\gamma f)}$.
\end{theorem}
We summarize Theorem \ref{boundedtheorem} by saying
 that the almost isomorphism identifies the thermodynamic
formalisms of the bounded relatively regular functions on
$S$ and $T$.

We now turn to somewhat more general results, which combine to give
Theorem \ref{boundedtheorem} as a corollary:
Propositions \ref{uniqueprop},
\ref{existprop} and \ref{propcorresp} imply Theorem
\ref{boundedtheorem}.

\section{Uniqueness}\label{sec:uniqueness}

In this section, we prove   the following uniqueness result.

\begin{prop}\label{uniqueprop}
Let $S$ be a Markov shift and $f\in \mathcal
E_1(S,\mathcal W) \cup \mathcal
E_{0+}(S,\mathcal W)$.
 Assume that $f$ is upper bounded and that the
pressure $P(S,f)$ is finite.

Then $(S,f)$ has at most one ergodic equilibrium measure $\mu$
giving positive measure to $\bigcup_{W\in\mathcal W} [W]$, and
this measure must have full support.

\end{prop}

We first note that restriction to a subsystem $S_{\mathcal W}$
does not affect the pressure.

\begin{lem} \label{samepressure}
Suppose $\mathcal W$ is a nonempty set of $S$-words and $f$ is a
Borel function from $S$ to $\mathbb R$. Let $f_{\mathcal W}$
denote the restriction of $f$ to $S_{\mathcal W}$. Then
$P(S,f)=P(S_{\mathcal W},f_{\mathcal W})$.
\end{lem}
\begin{proof}
The inequality $P(S,f)\geq P(S_{\mathcal W},f_{\mathcal W})$
 is trivial.
 For the
reverse inequality, let $\mu_n$ be
 a  sequence of ergodic measures
 of $S$ with $\lim_n P(S,f,\mu_n)=P(S,f)$ and choose
a word $W$ in $\mathcal W$. Approximate $\mu_n$ by an ergodic
measure $\nu_n$
 such that $P(S, f, \nu_n) \geq P(S,f,\mu_n )-1/n$,  and  also
 $\nu_n([W])>0$. (For example, for $\nu_n$ use a $k$-step Markov
perturbation of the $k$-step Markov approximation to $\mu_n$,
for sufficiently large $k$.)
Then $P(S_{\mathcal W},f_{\mathcal W})\geq \lim_n P(S, f,
\nu_n)=P(S,f)$.
\end{proof}

\begin{proof}[Proof of \ref{uniqueprop}]
 Let $f\in\mathcal E_1(S,\mathcal W)\cup
\mathcal E_{0+}(S,\mathcal W)$, for some collection
of words $\mathcal W$. Assume that $f$ is upper bounded and that
$P(S,f)<\infty$.

To  prove uniqueness, we assume that $\mu_1$ and $\mu_2$ are two
ergodic equilibrium measures assigning positive measure to
$\cup_{W\in \mathcal W}[W]$,
 and show that
they must coincide. Choose  words $w_1,w_2\in\mathcal W$ such that
$\mu_1([w_1])>0$ and $\mu_2([w_2])>0$. We assume that these words
have the same length $L$,  by lengthening the shorter one if
necessary. For $i=1,2$, we choose nonempty words $a_i,b_i$ such
that $\mu_i([W_i])>0$ for $W_i=w_ia_iw_ib_i$.
 We arrange it so that $W_1$ and $W_2$ have a common
length $N$.

Let $\mathcal G$ be the infinite graph whose vertices are the
$S$-words of length $N$, with an edge from vertex $u_0\cdots
u_{N-1}$ to vertex $u_1\cdots u_N$ iff $u_0u_1\cdots u_N$ is an
$S$-word, and we label such an edge  $u_0$. As is well known, this
graph determines a Markov shift, and the labeling defines a
block code to $S$ which is a topological conjugacy, i.e., a
shift-commuting homeomorphism.

Now define another graph $\overline{\mathcal G}$ as follows.
 $\overline{\mathcal G}$ will include two distinguished
vertices    $\overline{v_1}$ and $\overline{v_2}$. For each path
$p=p_1\cdots p_k$ in $\mathcal G$ with  initial vertex $W_i$ and
terminal vertex $W_j$, and  no intermediate vertices equal to
$W_i$ or $W_j$, $\overline{\mathcal G}$ contains a path $\overline
p$ of equal length from
 $\overline{v_i}$ to $\overline{v_j}$. These paths contain
all edges of  $\overline{\mathcal G}$ and overlap only at their
initial and terminal vertices in
 $\{ \overline{v_1},\overline{v_2}\}$. The $S$-word of length
$k$ labeling the edges $p_1,p_2, \dots , p_k$ is used to likewise
label the $k$ edges
 $\overline{p}_1,\overline{p}_2, \dots , \overline{p}_k$.
Let $R$ be the Markov shift defined from the resulting graph and
let $\varphi :R\to S$ be the injective block code defined by
the edge labeling. The image of $\varphi$  has measure one for
both $\mu_1$ and $\mu_2$. Let $\overline f$ denote the function
 $f\circ S^L\circ \varphi$.

 Let $\overline{\mathcal A}$ denote the alphabet of $R$. We claim
that
\begin{equation}\label{relclaim}
   \overline f\in\mathcal E_1(R,\overline{\mathcal A})
\cup \mathcal E_{0+}(R,\overline{\mathcal A})
\end{equation}
(contrarily to $f$ on $S$). Suppose $\overline x, \overline y$ are
sequences in $R$ with $\overline x[-m,n]= \overline y[-m,n]$ for
some integers $m,n\geq0$ which we may and do assume to be taken
maximum (the case where
 $\overline x$ and $\overline y$ agree in all nonnegative coordinates
or all nonpositive coordinates is similar and is left
to the reader). The initial vertex of $\overline x_{-m}$ and the
terminal vertex of $\overline x_n$ are contained in $\{ \overline
v_1, \overline v_2\}$. Let $x=\varphi (\overline x)$ and
 $y=\varphi (\overline y)$. Then
 for some $i,j\in\{1,2\}$,
$$   x[-m,-m+L-1]=y[-m,-m+L-1]=w_i$$
 and since $x[n,n+N-1]=w_ja_jw_jb_j$,
 $$
  x[n+N-|w_jb_j|,n+N-|b_j|-1]=
  y[n+N-|w_jb_j|,n+N-|b_j|-1]=w_j.
 $$
Remark that
\begin{align*}
   (-m+L-1)-L \ &\ \leq \ -1-m \ , \ \ \text{ and }\\
  (n+N-|w_jb_j|)-L \ &\ =\ n+|a_j|\ \geq \ n+1
\end{align*}
so $f\in\mathcal E_{p}(S,\{w_1,w_2\})$ gives:


\[
 |\overline
f(\overline x) - \overline f(\overline y)|\ =\ |f(S^Lx) - f(S^Ly)|
\ \leq \ \omega_{m+L} + \omega_{n+M}
\]
where  $M=N-|b_j|-L-1\geq 0$, and the sequence $(\omega_n )$ comes
from Definition \ref{psummdef} for $f$. Now define $\overline
\omega_n = \max \{\omega_{n+L},\omega_{n+M}\}$, so
\[
 |\overline
f(\overline x) - \overline f(\overline y)|\ \ \leq \
 \overline \omega_m + \overline \omega_n \
\]
and
\[
\sum_{n=1}^{\infty}  n^p \overline \omega_n \ \leq \
\sum_{n=1}^{\infty} (n+L)^p\omega_{n+L} + (n+M)^p \omega_{n+M} \
<\  \infty \ .
\]
Moreover if $f$ depends only on future coordinates, then this is
obviously the case for $\overline f$. Therefore  the claim
(\ref{relclaim}) is proved (in particular $\overline f$ has
bounded oscillation on cylinders $[a]$,
$a$ a vertex of $\overline{\mathcal
G}$, contrarily to $f$).

Next, in the case where $f$ also depends on past coordinates, we
make the following  standard replacement to obtain a function
depending only on future coordinates. Following \cite{BowenLNM},
let $h:R\to\mathbb R$ be given by :
 $$
    h(x) = \sum_{k\geq0} \overline f(R^kx)-
\overline f(R^kx^-)\qquad \forall x\in R
 $$
where $x^-\in R$ with $x^-_i=x_i$ for $i\geq0$ and
$(x^-_i)_{i\leq0}$ depending only on $x_0$. The above series is
bounded by $\sum_{k\geq0} \omega_{k}$, and converges. One easily
checks that $h$ is bounded and has summable variations. One also
sees that $\overline f+h\circ R-h$ depends only on
$(x_i)_{i\geq0}$.

If $f$ depends only on future coordinates, let $h=0$. Then in both
cases $\overline f+h\circ R-h$ is upper bounded, has summable
variations and depends only on future coordinates. We can apply
the main theorem of \cite{BS} and see that $(R,\overline f)$ has
at most one ergodic equilibrium measure, and that this measure has
full support when it exists.

To bring back this result to  $(S,f)$ first notice that
$\mu(\overline f+h\circ R-h)=\mu(f\circ\varphi)$ for all
$R$-invariant probability measures. We claim further that
$P(R,\overline f)=P(S,f)$.
 This follows from Lemma
\ref{samepressure} because $\varphi$ is a continuous bijection,
with Borel measurable inverse function, and the image of $\varphi$
is  $S_{\{W_1,W_2\}}$. Thus $\mu_1$ and $\mu_2$ must correspond
under $\varphi$ to the unique equilibrium measure of $(R,\overline
f)$. Hence $\mu_1=\mu_2$.

Because $\varphi$ is continuous and has dense image, it also
follows that $\mu_1=\mu_2$ has full support in $S$.
\end{proof}

\section{Existence}\label{sec:existence}

The issue of existence in our context is more complicated. As in
\cite{Sarig0,Walters0,Yuri0}, we will use the concept of a
{\sl weak equilibrium measure} for $(S,f)$:
 an invariant probability measure $\mu$
for which $P(S,f,\mu )$  is not necessarily defined, but which
for some measurable function $h$ satisfies the following:
\[
  -\sum_{a\in\mathcal A} 1_{[a]} \log E_\mu(1_{[a]}|
    S^{-1}\mathcal B) + f+h\circ S-h
                  \in L^1(\mu)
\]
and
 \begin{equation} \label{eq-generalizedP}
    \int \left( -\sum_{a\in\mathcal A} 1_{[a]} \log E_\mu(1_{[a]}|
    S^{-1}\mathcal B) + f+h\circ S-h \right) \, d\mu = P(S,f),
 \end{equation}
where $\mathcal A$ is the set of states of $S$, $\mathcal B$ is
the $\sigma$-algebra of Borel measurable subsets,
$E_\mu(\cdot|\cdot)$ is the conditional expectation.
(The function $h$ was
assumed to be locally H\"older-continuous in \cite{Sarig0}, but
not here.)

\begin{remark}\label{weakeq}
It is possible for a positive recurrent H\"older continuous
potential to have a weak equilibrium measure when there is no
equilibrium measure
 \cite[Sec.7]{Sarig0}. When $S$ has finite topological entropy and
$f$ is bounded, then a weak equilibrium measure for $(S,f)$ must
be an equilibrium measure. However there exist upper-bounded
potentials which are positive recurrent, have finite pressure and
summable variations, define a weak equilibrium measure which has
finite entropy, and, yet, this measure is not a (strong)
equilibrium measure.

 O. Sarig pointed out to us the following example of
such a potential. Fix $1/2<\alpha<1$ and take the interval map
$T:[0,1]\to[0,1]$ defined by $T(x)=x(1+2^\alpha x^\alpha)$ for
$x<1/2$ and $T(x)=2x-1$ otherwise \cite{LSV}. Let
$\varphi(x)=-\log|T'|-\alpha\log (x/Tx)\leq 0$. Using the
partition of $[0,1]\setminus \bigcup_{k\geq0}
(T^{-k}(1/2)\cap[0,1/2))$ into its connected components, one
obtains from $(T,\varphi)$ a Markov shift with a potential which
gives the required example (use \cite{Sarig4,Sarig3} for the
summable variations and \cite{Thaler} for estimating the density
of the weak equilibrium measure). We do not know if one can find
an example where additionally the potential is
H\"older-continuous.
\end{remark}

We now state our main existence result.

\begin{prop}\label{existprop}
Let $S$ be a Markov shift and $f\in \mathcal
E_1(S,\mathcal W) \cup \mathcal
E_{0+}(S,\mathcal W)$.
 Assume that $f$ is upper bounded and the
pressure $P(S,f)$ is finite.

If $(S,f)$ is positive recurrent, then there exists a weak
equilibrium measure. Conversely, if there exists an equilibrium
measure $\mu$
 with $\mu S_{\mathcal W}=1$, then $(S,f)$
is positive recurrent.

If $f$ is bounded  and $\mu S_{\mathcal W}=1$ whenever
$\mu$ is an equilibrium measure for $f$,  then $(S,f)$ is positive
recurrent if and only if $f$ has an equilibrium measure.
\end{prop}

\begin{remark}
Propositions \ref{uniqueprop} and \ref{existprop}
are known (in the case $S_{\mathcal W}=S$)
if $f$ is uniformly
locally constant \cite{GS} or if it is uniformly H\"older
continuous or more generally has summable variations and depends only
on the future \cite{Sarig0}.
We note that Gibbs measures (see \cite{MU1}) are a rather different
issue \cite{Sarig3}.
\end{remark}

\begin{remark}
 We do not know if the existence of a weak equilibrium
measure implies that a potential in $\mathcal E_1(S,\mathcal W)
\cup \mathcal E_{0+}(S,\mathcal W)$ is positive recurrent.
\end{remark}

\begin{remark}
The paper \cite{FFY} gives
more general sufficient conditions
 for existence of equilibrium measures
of continuous functions on a Markov shift.
We have not exploited those conditions here because
we do not know when the equilibrium measures produced in \cite{FFY}
have full support, and we do not know whether
the  $Z$-recurrence condition in \cite{FFY} must be
preserved under passage to a smaller system $S_W$.
\end{remark}

\begin{remark}
 $f\in\mathcal E_p(S,\mathcal W)$ is
arbitrary on the subshift $S_*$ of $S$ obtained by excluding all
words in $\mathcal W$. Hence, the  exclusion
in Definition \ref{goodclass}
of equilibrium
measures living on $S_*$ is necessary
for drawing conclusions about equilibrium measures on $S$.
Also, in the definition we could equivalently
 require $f$ to be positive recurrent on  $S_{\mathcal W}$
 rather than
on $S$, since given $f$ on $S_{\mathcal W}$
we could extend $f$ to $S$ without changing the pressure or
set of equilibrium measures
\end{remark}

\begin{proof}[Proof of \ref{existprop}]
We continue with the notation and constructions used in the
proof above for Prop. \ref{uniqueprop}. Let
$(R,\overline f)$ be constructed as in that proof,
using as $w_1=w_2$ some word from $\mathcal W$ to be specified.

Assume first that $(S,f)$ is positive recurrent. Take
$w_1=w_2\in\mathcal W$ arbitrarily. Positive recurrence of
$(R,\overline f)$ follows from $P(R,\overline f)=P(S,f)$ and
the coincidence of local partition functions for all $n\geq1$,
 $$
    Z_n(S,f,W_1)=Z_n(R,\overline f,\overline v_1)
 $$
 where $W_1$ is a word in $\mathcal W$ and
$\overline{v}_1$ is an element of the alphabet of $R$, both
defined as in the proof for Prop. \ref{uniqueprop}.
 Then $(R,\overline f)$  has a weak equilibrium measure
by Theorem 7 of \cite{Sarig0}. $\overline\varphi$ being
one-to-one, this gives a weak equilibrium measure for $(S,f)$.

Conversely, assume now that $(S,f)$ has an equilibrium measure
$\mu$  with $\mu(S_W)=1$. By assumption, we can choose
some word $w_1=w_2\in\mathcal W$ with $\mu([w_1])>0$. $\mu$ can
then be pulled back to $R$ by $\overline\varphi$. As
$P(R,\overline f)= P(S,f)$  the pullback of $\mu$ is an
equilibrium measure for $(R,\overline f)$. This implies the
positive recurrence of $(R,\overline f)$ according to
\cite{Sarig0} and \cite{BS} (\cite{BS} says that any equilibrium
measure satisfies Sarig's Ruelle-Perron-Frobenius theorem;
\cite{Sarig0} says that existence of such an invariant measure
implies positive recurrence). The positive recurrence of $(S,f)$
follows as above.

It remains to show that when $f$ is bounded, a weak equilibrium
measure for $f$ must be an equilibrium measure. This is an
immediate consequence of
Lemma \ref{lemWeakEqM} below.
\end{proof}

\begin{lem} \label{lemWeakEqM}
Let $f\in\mathcal E_0(S,\mathcal W)$ be upper-bounded with
$P(S,f)<\infty$. If $\mu$ is a weak equilibrium measure for
$(S,f)$ such that $f\in L^1(\mu)$ and
$\mu\left(\bigcup_{W\in\mathcal W} [W]\right)>0$ then $\mu$ is an
equilibrium measure.
\end{lem}

\begin{proof}[Proof of Lemma \ref{lemWeakEqM}]
Let $h$ be the measurable function given by the assumption that $\mu$
is a weak equilibrium measure. Recall our notation
$S_nf=f+f\circ S+\dots+f\circ S^{n-1}$.

The first point to see is that $\int h-h\circ S\, d\mu=0$. We follow
Ledrappier \cite{Ledrappier} by observing that $h-h\circ S$ is
integrable as the difference of two integrable functions, namely
$f$ and $f+(h-h\circ S)$, and
therefore Birkhoff's ergodic theorem gives
 $$
    \int h-h\circ S \, d\mu = \lim_{n\to\infty} \frac1nS_n(h-h\circ S)
     = \lim_{n\to\infty} \frac1n(h-h\circ S^n)\ \ \
\text{a.e.}
 $$
Recurrence
implies that $0$ is an accumulation point for the sequence
$\frac1n(h-h\circ S^n)$ a.e. Therefore the above limit is
zero, proving the first point and therefore
 \begin{equation}\label{eq-A}
    \int f\, d\mu = \int f+h-h\circ S\, d\mu \ .
 \end{equation}

We have the following properties:
 \begin{itemize}
  \item $P(S,f)<\infty$;
  \item $f\in L^1(\mu)$;
  \item because of the eventual summable variation property,
there exists a cylinder $[W_*]$ with $\mu([W_*])>0$ such that
for all $x,y\in[W_*]$, $n\geq0$ with $S^n(x),S^n(y)\in[W_*]$,
$|S_nf(x)-S_nf(y)|<\text{const}$.
 \end{itemize}
The last property is proved as in the proof of Proposition
\ref{uniqueprop} (as there, we replace $f$ by $f\circ S^L$). We
claim that
 \begin{equation}\label{eq-B}
    \int -\sum_{a\in\mathcal A} 1_{[a]} \log E_\mu(1_{[a]}|
    T^{-1}\mathcal B)  \, d\mu \ =\ h(S,\mu)
 \end{equation}
despite the fact that the partition $\{[a]: a\in\mathcal A\}$ may
have infinite entropy. This is proved in
\cite[p. 1389]{BS},
assuming the three
properties listed above.
Because $\mu$ is a weak equilibrium measure,
the facts (\ref{eq-A}) and (\ref{eq-B})
together show that $\mu$ is indeed an equilibrium measure.
\end{proof}


\section{Correspondence}\label{sec:correspondence}

\begin{definition}
 $ \mathcal E^{u\flat}(S)$ is the class of
upper bounded functions $f$ with
$P(S,f)<\infty$ such that
 there exists
$\mathcal W$ (allowed to depend on $f$) such that $ f\in \mathcal
E_1(S,\mathcal W)\cup \mathcal
E_{0+}(S,\mathcal W)$ and for any equilibrium measure
$\mu$ of $f$, $\mu [W]>0$ for some $W\in \mathcal W$.
\end{definition}

\begin{proposition}\label{propcorresp}
 Suppose $R\to S$, $R \to T$ is an almost isomorphism of
Markov shifts $S$ and $T$. This almost isomorphism induces
 mod-s.e. bijections
 $\ \Gamma: \mathcal E_p(S) \to
\mathcal E_p(T) $ and
 $\ \Gamma: \mathcal E^{u\flat}(S) \to
\mathcal E^{u\flat}(T) $.
Moreover, for $f\in  \mathcal E^{u\flat}(S)$ we have:
 \begin{itemize}
  \item $(S,f)$ is positive recurrent if and only if $(T,\Gamma(f))$
is positive recurrent;
  \item an invariant probability measure $\mu$ of
$S$ is an equilibrium measure, resp. fully supported weak
equilibrium measure, for $(S,f)$ if and only if $\gamma \mu$ is
an equilibrium measure, resp. fully supported weak equilibrium
measure, for $(T,\Gamma(f))$.
\end{itemize}
\end{proposition}

\begin{remark}
We do not know whether a weak equilibrium measure for a H\"older
continuous potential is always fully supported.
\end{remark}

\begin{proof}
Let $\gamma:K\subset S\to K'\subset T$ be the bijection given by
the almost isomorphism according to Proposition \ref{magicprop}.
We first prove that if $f\in \mathcal E_p(S,\mathcal W)$
then for some collection $\mathcal V$ of $T$-words,
 \begin{equation}\label{eq-gErel}
   g=f\circ\gamma^{-1}:K'\to\mathbb R \in\mathcal
E_p(T,\mathcal V)\ .
 \end{equation}
(If one likes, one can extend $g$ to the whole of $T$.)

Let $W_0\in\mathcal W$. Let $V_1$ be some magic $T$-word for
the map $R\to T$ such that $K'$ contains all points in which
$V_1$ occurs infinitely often in the future and the past.
Choose a $T$-word $V_2$ and an integer
$J\geq0$ such that for all $y\in K'$ and $i\in \mathbb Z$,
\begin{align*}
\text{if }\  &y[i,i+|V_1V_2V_1|-1] = V_1V_2V_1 \ , \\
 \text{ then }\
&(\gamma^{-1}y)[i+J,i+J+|W_0|-1]=W_0\ .
\end{align*}
 We may also assume $J\leq |V_1V_2V_1|-|W_0|$.
Let $V$ be the $T$-word  $V_1V_2V_1$. We claim that $g\in \mathcal
E_p(T,\mathcal V)$ with $\mathcal V=\{V\}$.

Indeed let $x,y\in K'$ and the integers $m,n$ be such that
$x[-m,n]=y[-m,n]$ begins and ends with
the word  $V$. Then
$x=\gamma(u)$ and $y=\gamma(v)$, where
\[
u[-m+J,n-|V|+J+|W_0|+1]\ =\ v[-m+J,n-|V|+J+|W_0|+1]
\]
and this word
begins and ends with $W_0$.
Therefore
\[
|g(x)-g(y)|\ =\ |f(u)-f(v)|\ \leq\
\omega_{m-J}
+\omega_{n-|V|+J}
\]
which implies
 (\ref{eq-gErel}) by an argument similar to the proof of
(\ref{relclaim}).

We next check that $\Gamma$ is well-defined on the level of
somewhere equivalent classes. Let $f_1,f_2\in\mathcal
E_p(S)$ be somewhere equivalent, i.e., for some word
$w$, $f_1=f_2$ over $S_w$. By enlarging the word $W_0$ in the
above construction to a word $W$ which contains  $w$, we obtain
that $T_{V}\subset\gamma(S_W)$. This implies that
$\Gamma(f_1)=\Gamma(f_2)$ s.e. Because  $\gamma:K\to K'$ is a
bijection, the induced map $\Gamma$ on  somewhere equivalent
classes is also a bijection.

We now check that for $f\in \mathcal E^{u\flat}(S)$, there
exists $g\in \mathcal E^{u\flat}(T)$ somewhere
equivalent to $f\circ\gamma^{-1}$.
We define $g$ to be
$f\circ\gamma^{-1}$ on $K'$.
Let $T_*$  denote the restriction of
$T$  to the complement $K'_*$
of $K'$; choose an upper bounded
 function $g_*$ on  $K'_*$ such that
$P(T_*,g_*)<P(T,g)$; and on  $K'_*$, define $g=g_*$. Now  any
equilibrium measure for $(g,T)$ must be supported on $K'$.
 Using Lemma \ref{samepressure}, we have
\[P(T,g)=P(T|{K'},g)=P(S|K,f|K)=P(S,f)<\infty\ .
\]
 This proves that $\Gamma(\mathcal
E^{u\flat}(S))\subset \mathcal
E^{u\flat}(T)$ modulo s.e. It follows that
$\Gamma:\mathcal E^{u\flat}(S)\to \mathcal
E^{u\flat}(T)$ is a bijection modulo s.e.

%

Note that Proposition \ref{uniqueprop} says that if $\mu$ is an
equilibrium measure for $f\in\mathcal E^{u\flat}(S)$, it is
ergodic and has full support. By Proposition \ref{magicprop},
$(T,\gamma\mu)$ is isomorphic to $(S,\mu)$ as soon as $\mu$ is
fully supported. Thus we see that the equilibrium measure $\mu$ of
$(S,f)$ corresponds to a measure $\gamma\mu$ with the same
pressure for $(T,\Gamma(f))$. Because $P(S,f)=P(T,\Gamma(f))$,
 $\gamma\mu$ is the equilibrium measure for
$(T,\Gamma(f))$ as claimed.

The same reasoning applies to fully supported weak equilibrium
measures.
\end{proof}

\end{document}